\input amstex
\input amsppt.sty   
\nologo
\hsize 30pc
\vsize 47pc
\magnification=\magstep1
\def\nmb#1#2{#2}         
\def\cit#1#2{\ifx#1!\cite{#2}\else#2\fi} 
\def\totoc{}             
\def\idx{}               
\def\ign#1{}             
\redefine\o{\circ}
\define\X{\frak X}
\define\al{\alpha}

\define\ga{\gamma}
\define\de{\delta}
\define\ep{\varepsilon}
\define\ze{\zeta}

\define\th{\theta}

\define\ka{\kappa}
\define\la{\lambda}
\define\rh{\rho}
\define\si{\sigma}

\define\ph{\varphi}

\define\om{\omega}
\define\Ga{\Gamma}

\define\La{\Lambda}
\define\Si{\Sigma}

\redefine\i{^{-1}}
\define\x{\times}
\define\g{\frak g}
\define\Id{\operatorname{Id}}

\define\ad{\operatorname{ad}}
\define\Der{\operatorname{Der}}

\define\Nor{\operatorname{Nor}}
\define\Fl{\operatorname{Fl}}
\let\on=\operatorname
\def\today{\ifcase\month\or
 January\or February\or March\or April\or May\or June\or
 July\or August\or September\or October\or November\or December\fi
 \space\number\day, \number\year}
\topmatter
\title  The Riemannian geometry of orbit spaces.\\
The metric, geodesics, and integrable systems
\endtitle
\author Dmitri Alekseevsky, Andreas Kriegl, 
Mark Losik, and Peter W. Michor  
\endauthor
\date February 14, 2001 \enddate
\rightheadtext{The Riemannian geometry of orbit spaces}
\leftheadtext{D\.V\. Alekseevsky, A\. Kriegl, M\. Losik, P\.W\. Michor}
\affil
Erwin Schr\"odinger International Institute of Mathematical Physics, 
Boltzmanngasse 9, A-1090 Wien, Austria
\endaffil
\address 
D.V. Alekseevsky: 
Department of Mathematics,
University of Hull,
Cottingham Road,
Hull, HU6 7RX,
England
\endaddress
\email
d.v.alekseevsky\@maths.hull.ac.uk 
\endemail
\address
A\. Kriegl: Institut f\"ur Mathematik, Universit\"at Wien,
Strudlhofgasse 4, A-1090 Wien, Austria
\endaddress
\email Andreas.Kriegl\@univie.ac.at \endemail
\address
P\. W\. Michor: Institut f\"ur Mathematik, Universit\"at Wien,
Strudlhofgasse 4, A-1090 Wien, Austria; {\it and:} 
Erwin Schr\"odinger Institut f\"ur Mathematische Physik,
Boltzmanngasse 9, A-1090 Wien, Austria
\endaddress
\email Peter.Michor\@esi.ac.at \endemail
\address
M. Losik: Saratov State University, ul. Astrakhanskaya, 83,
410026 Saratov, Russia
\endaddress
\email LosikMV\@info.sgu.ru \endemail
\thanks  
P.W.M. was supported  
by `Fonds zur F\"orderung der wissenschaftlichen  
Forschung, Projekt P~14195~MAT'.
\endthanks 
\keywords orbit spaces, geodesics, invariants, representations \endkeywords
\subjclass 26C10, 53C22 \endsubjclass
\abstract 
We investigate the rudiments of Riemannian geometry on orbit spaces 
$M/G$ for isometric proper actions of Lie groups on Riemannian 
manifolds. 
Minimal geodesic arcs are length minimising curves in the metric 
space $M/G$ and they can hit strata which are more singular only at the end 
points. This is phrased as convexity result. The geodesic spray, 
viewed as a (strata-preserving) vector field on $TM/G$, leads to the 
notion of geodesics in $M/G$ which are projections under $M\to M/G$ 
of geodesics which are normal to the orbits. It also leads to 
`ballistic curves' which are projections of the other geodesics. 
In examples (Hermitian and symmetric matrices, and more generally polar 
representations) we compute their equations by singular symplectic 
reductions and obtain generalizations of Calogero-Moser systems with 
spin. 
\endabstract
\endtopmatter

\document

\head Table of contents \endhead

\noindent 1. Introduction \leaders 
\hbox to 1em{\hss .\hss }\hfill {\eightrm 1}\par   
\noindent 2. Stratification of orbit spaces 
\leaders \hbox to 1em{\hss .\hss }\hfill {\eightrm 3}\par 
\noindent 3. The orbit space $M/G$ as a metric space 
\leaders \hbox to 1em{\hss .\hss }\hfill {\eightrm 4}\par 
\noindent 4. Vector fields, geodesics, and Jacobi fields on orbit spaces 
\leaders \hbox to 1em{\hss .\hss }\hfill {\eightrm 6}\par 
\noindent 5. Example: Hermitian and symmetric matrices 
\leaders \hbox to 1em{\hss .\hss }\hfill {\eightrm 9}\par 
\noindent 6. Ballistic curves on polar representations 
\leaders \hbox to 1em{\hss .\hss }\hfill {\eightrm 15}\par 


\head\totoc\nmb0{1}. Introduction \endhead

Differential geometry deals ordinarily with smooth objects on smooth 
manifolds. However, in many cases this assumption is too restrictive. 
In various contexts different types of non-smooth manifolds or 
`manifolds with singularities' appear naturally, like complex 
algebraic varieties, orbifolds, $V$-manifolds \cit!{24}, \cit!{7}, 
limits of Riemannian manifolds with respect to the Gromov-Hausdorff 
metric \cit!{9}, etc.
In the present paper we consider a rather special but interesting and 
important class of non-smooth manifolds, namely the orbit spaces 
$M/G$ of Riemannian manifolds $M$ with respect to groups $G$ of 
isometries. Our goal is to develop the Riemannian geometry of orbit 
spaces. An orbit space has the structure of a stratified manifold 
with smooth strata, namely the connected components of the sets 
$(M/G)_{(H)}$ of orbits of some given orbit type (or isotropy type) 
$(H)$. 

Let us describe the structure of the paper. The basic facts about the 
orbit stratification of an orbit space $M/G$ is presented in section 
\nmb!{2}. 

An orbit space of a complete Riemannian $G$-manifold $M$ has 
a natural structure of a metric space. Some elementary properties of 
this metric structure are collected in section \nmb!{3}. In 
particular, we establish some convexity properties of the 
stratification with respect to distance minimizing curves which we 
call minimal geodesic arcs. For 
example, the set $(M/G)_{\text{reg}}$ of regular  
orbits is a convex open dense submanifold of $M/G$.

Some basic global differential geometric objects on $M/G$ are defined 
in section \nmb!{4}. We define the algebra of smooth functions 
$C^\infty(M/G)$ on $M/G$ as the algbera $C^\infty(M)^G$ of 
$G$-invariant smooth functions on $M$, and the Lie algebra of vector 
fields $\X(M/G)$ as the Lie algebra of all derivations of 
$C^\infty(M/G)=C^\infty(M)^G$ which preserve all ideals of functions 
vanishing on some stratum. The deep result \cit!{26} of G\.~Schwarz 
states that the natural homomorphism 
$\X(M)^G\to \on{Der}(C^\infty(M/G))$ has $\X(M/G)$ as image. 
Each element of $\X(M/G)$ induces a derivation on the algebra of smooth 
functions on each stratum and thus defines an 
ordinary vector fields along each stratum. In order to define geodesics 
which connect different strata we consider the geodesic spray as a vector 
field on $TM/G$, the projection of the geodesic spray 
$\Ga\in \X(TM)^G$. The integral curves of the geodesic spray $\Ga$ on 
$TM/G$ stay within the strata, but their projection to $M/G$ connects 
different strata in $M/G$. We define geodesics in $M/G$ as 
projections to $M/G$ of integral curves of $\Ga$ on $TM/G$ which are 
orthogonal to orbits, i.e., in $TM$ their initial vector should be in 
in $\Nor(M)_x=T_x(G.x)^\bot$ for $x\in M$.
More generally, the projection on $M/G$ of an arbitrary integral 
curve of the geodesic spary $\Ga$ on $TM/G$ is called a `ballistic 
curve'. Imitating the classical Riemannian case, we establish some 
properties of the geodesic spray on $TM/G$ and define Jacobi fields 
along gedesics and along ballistic curves. 

Section \nmb!{5} is devoted to a more systematic study of ballistic 
curves in $M/G$ for the simple model case when the unitary group 
$G=SU(n)$ acts on the space $M=H(n)$ of Hermitian matrices by 
conjugation. Using the singular Hamiltonian reduction by Sjamaar and 
Lerman \cit!{27}, we derive the Hamiltonian equation for a ballistic 
curve in $H(n)/SU(n)$. It turns out that it is a Calogero-Moser 
system with spin. In the special case when the momentum 
$Y\in \frak s\frak u(n)^*$ has maximal isotropy group, i.e.\ when 
$Y=\sqrt{-1}(c.1_n+w\otimes w^*)$ where $w$ is a vector with 
$|w|^2=-c>0$ the equation is the classical Calogero Moser system, 
which reproduces results of Kazhdan, Kostant, and Sternberg, 
\cit!{10}.

The last section \nmb!{6} is devoted to the generalization of the  
results on ballistic curves to the orbit space $V/G$ where $V$ is an 
Euclidean space and $G\subset SO(V)$ is a connected subgroup of the 
orthogonal group whose action on $V$ is polar, i.e.\ admits a 
section.  

We thank Nikolai Reshetikhin for helpful discussions and pointing out 
references \cit!{10}, \cit!{3}, and \cit!{11}.

\head\totoc\nmb0{2}. Stratification of orbit spaces \endhead

\subhead\nmb.{2.1}. The setup \endsubhead
Let $M$ be a connected $G$-manifold, where $G$ is a Lie group. 
The $G$-action is called \idx{\it proper} if $G\x M\to M\x M$, given by 
$(g,x)\mapsto (g.x,x)$, is a proper mapping. 
It is well known that proper $G$-actions admit \idx{\it slices} 
(\cit!{20}, \cit!{21}, \cit!{16}): For 
each $x\in M$ there exists a submanifold $S_x\subset M$ containing 
$x$, an open $G$-invariant neighborhood $W$ of $G.x$, and a smooth 
equivariant retraction $r:W\to G.x$ such that $S_x=r\i(x)$. 
Moreover, for the isotropy group $G_x$ of $x$ we have 
$G_x.S_x\subseteq S_x$, and $g\in G$ with 
$g.S_x\cap S_x\ne \emptyset$ must lie in $G_x$. Moreover, the slice 
$S_x$ is a manifold and $G_s\subseteq G_x$ for each $s\in S$. 
Finally, the action $G\x S_x\to W$ induces an $G$-equivariant 
diffeomorphism $G\x_{G_x}S_x\to W$. This implies that 
$C^\infty(W)^G=C^\infty(S_x)^{G_x}$ via restriction.

We consider the orbit space $M/G$ and the canonical projection 
$\pi:M\to M/G$, and we endow the quotient space $M/G$ with the 
following smooth structure:
The quotient topology and the sheaf of smooth real valued functions  
$U\mapsto C^\infty(U):=C^\infty(\pi\i(U))^G$. A mapping 
$\ph:M/G\to M'/G'$ is called \idx{\it smooth} if it respects these 
sheafs. For a slice $S_x$ as above  on the orbit spaces we have 
$C^\infty(W/G)=C^\infty(S_x/G_x)$. Therefore the local smooth 
structure of $M/G$ coincides with the smooth structure of $S_x/G_x$. 
A mapping $\ph:M/G\to M'/G'$ is smooth if and only if 
$\ph^*C^\infty(M'/G')\subseteq C^\infty(M/G)$.

\subhead\nmb.{2.2} \endsubhead
Let $H$ be a closed subgroup of $G$ and let $(H)$ denote the 
conjugacy class of $H$. For two closed subgroups $H_1$ and $H_2$ we 
write $(H_1)\le(H_2)$ if $H_1$ is conjugated to a subgroup of $H_2$. 

Let $M_{(H)}$ denote the set of points of $M$ whose isotropy groups 
belong to $(H)$. It is known that $M_{(H)}$ is a smooth submanifold 
of $M$ for proper actions (see \cit!{16}, 7.4). Put 
$(M/G)_{(H)}:=\pi(M_{(H)})=(M_{(H)})/G$, call this the isotropy  
stratum of type $(H)$, and call any connected component of 
this an \idx{\it orbit stratum} of $M/G$.

\proclaim{Proposition} \cit!{26}, \cit!{16}, \cit!{21}.
\roster
\item The isotropy stratum $(M/G)_{(H)}=M_{(H)}/G$ 
       is a smooth manifold, 
       the inclusion $(M/G)_{(H)}\to M/G$ is smooth, and 
       $\pi:M_{(H)}\to(M/G)_{(H)}$ is a smooth fiber bundle with fiber 
       type $G/H$.
\item We have a smooth fiber bundle $M_{(H)}\to G/N_G(H)$ where 
       $N_G(H)$ is the normalizer of $H$ in $G$, and where the 
       fiber over $g.N_G(H)$ is the fixed point set 
       $M^{gHg\i}\cap M_{(H)}$. See \cit!{16}, 7.3.  
\item The orbit strata of $M/G$ 
       form a locally finite partition of $M/G$.
\endroster
\endproclaim

\proclaim{\nmb.{2.3}. Theorem} 
\cit!{16}, 6.15ff.

There exists a unique minimal isotropy type $(K)$ such that 
\roster
\item $(M/G)_{(K)}$ is connected, locally connected, open, and dense 
       in $M/G$. 
\item The slice representations at all points of $M_{(K)}$ are 
       trivial, i.e.\ $G_x$ acts trivially on $S_x$ for all 
       $x\in M_{(K)}$.
\item $\dim(M/G)_{(K)} = \dim M - \dim G + \dim K$. 
\endroster
\endproclaim
The stratum $(M/G)_{(K)}$ is called the principal isotropy 
stratum and is denoted by $(M/G)_{\text{reg}}$. Likewise we write 
$M_{\text{reg}}:=M_{(K)}$. 

\subhead\nmb.{2.4} \endsubhead
Let $G$ be a compact group and let $\rh:G\to GL(V)$ be an orthogonal 
representation of $G$ on a real finite dimensional Euclidean vector 
space $V$.  
Let $\si=(\si_1,\dots,\si_n):V\to \Bbb R^n$, where $\si_1,\dots,\si_n$ is 
a system of generators for the algebra $\Bbb R[V]^G$ of invariant 
polynomials on $V$. The mapping $\si$ is proper and induces a 
homeomorphism between $V/G$ and the closed subset 
$\si(V)\subset \Bbb R^n$, see \cit!{25}. Since $\si$ is a polynomial 
map, $\si(V)$ is a semialgebraic subset, see \cit!{23} for an 
explicit description by polynomial equations and inequalities. 
By \cit!{28} and \cit!{13} the semialgebraic set $\si(V)$ has a 
canonical stratification into smooth algebraic submanifolds, called 
the \idx{\it Whitney stratification}.
The strata for this stratification are the connected 
components of the images under $\si$ of the set of points in $V$ 
where the rank of the system of polynomials $\si_1,\dots,\si_n$ is 
constant. 

\proclaim{Theorem}
\roster
\item \cit!{5} The mapping $\si:V\to \si(V)$ induces a bijection 
       $\bar\si:V/G\to\si(V)$ which maps the components of the 
       isotropy strata of $V/G$ diffeomorphically onto the strata of 
       $\si(V)$ as a semialgebraic set.
\item \cit!{25} $\si^*:C^\infty(\Bbb R^n)\to C^\infty(V)^G$ is a 
       surjective homomorphism of algebras. There exists a 
       continuous linear map 
       $\ph:C^\infty(V)^G\to C^\infty(\Bbb R^n)$ with 
       $\si^*\o\ph=\Id_{C^\infty(\Bbb R^n)}$, \cit!{15}.
\endroster
\endproclaim

\subhead\nmb.{2.5}. Remark\endsubhead 
Let $M$ be a smooth proper $G$-manifold. 
Then the orbit stratification of $M/G$ is locally given as the 
Whitney stratification of a semialgebraic subset in a vector space; 
by \cit!{14} this is in turn determined by the algebra of smooth 
functions on it. Thus the orbit stratification of $M/G$ is 
determined by $C^\infty(M)^G=C^\infty(M/G)$.
 
\head\totoc\nmb0{3}. The orbit space $M/G$ as a metric space 
\endhead   

Let $(M,g)$ be a connected complete Riemannian manifold and let $G$ 
be a Lie  group of isometries  which acts properly on $M$ (or 
equivalently, is closed in the full  group of isometries). Then we 
say that $(M,g)$  is a complete Riemannian $G$-manifold. Denote by 
$\pi : M \to \bar M = M/G$  the natural projection of $M$ onto the  
orbit space $M/G$. 
 
Denote by $d$ the natural metric structure on $M/G$ induced by 
the Riemannian metric $g$ of $M$. By definition 
the distance $d(\bar p, \bar q)$  is the minimum of the lengths of all  
curves in $M $ which connect the orbits $\bar p, \bar q$. 
  
Recall that a metric space $(X,d)$ is called to be of 
\idx{\it inner type} \cit!{2} or a \idx{\it path metric space} 
\cit!{9}
if the distance between any two points $p,q$ is equal to the length  
of a curve $pq$ connecting these points. Such a curve is called a 
\idx{\it minimal geodesic segment}. 

\proclaim{\nmb.{3.1}. Proposition} 
Let  $(M,g)$  be a complete  
Riemannian $G$-manifold. Then the following holds:
\roster 
\item The orbit space $(M/G, d)$ with the natural metric is a 
      complete metric space and a path metric space. 
\item Any minimal geodesic segment of $M/G$ is the projection of a normal 
      (i.e\. orthogonal to orbits) geodesic segment of $M$ .  
\item For every $\bar p \in M/G$ there exists $r>0$ such that each 
      $\bar q$ with $d(\bar p,\bar q)<r$ can be connected to $\bar p$ 
      by a unique minimal geodesic segment.  
\endroster 
\endproclaim 

Note that even if $M$ is compact, in \therosteritem3 one cannot 
choose the same $r>0$ for all points $\bar p\in M/G$ in general, 
see example \nmb!{3.3}.
 
\demo {Proof} 
\therosteritem1 and \therosteritem2.
For $p,q\in M$ there is a point $g.q\in G.q$ such that 
$d(p,g.q)$ is the distance from $p$ to the (closed) orbit $G.q$. 
Since $M$ is a complete Riemannian manifold, there is a geodesic $c$ of 
minimal length from $p$ to $g.q$. This geodesic is orthogonal to the 
orbit $G.p$: otherwise, by Gauss' lemma, we could find a shorter 
broken geodesic from $p$ to $g.q$. 
By a well known lemma (\cit!{16},~8.1), $c'(t)$ is orthogonal to each 
orbit which it meets. But then the length of $\pi\o c$ in 
$M/G$ equals the length of $c$ in $M$, which is the distance 
between the orbits $G.p$ and $G.q$ and thus equals 
$d(\bar p,\bar q)$.

The metric space $(M/G,d)$ is complete since each bounded closed 
set is compact: Use that the image of a closed geodesic ball 
$B_r(x)\subset M$ of radius $r$ with center $x\in M$ is the closed 
ball $B_r(\pi(x))$ of the same radius in $M/G$. 

\therosteritem3
Let $\bar p, \bar q$ are two points of $M/G$ with sufficiently  
short distance $d =d(\bar p, \bar q)$ and let $\ga, \de$ be  
two geodesic segments in $M$ of length $d$ connecting the orbits 
$\bar p, \bar q$.  
Transforming one of the geodesics by an appropriate isometry, we may 
assume that the geodesics start from the same point $p \in \bar p$. 
Since the geodesics are normal to the orbit $\bar p$ , they belong 
to the slice $S = \operatorname{exp}_pN$, where $N$ is a  
neighbourhood of the origin in the normal space of the orbit $\bar p$ 
at $p$. The end points $\ga(d), \de(d)$  of these geodesics belong 
to the orbit $\bar q$. Hence, by the main property of a slice, there 
exist an isometry $h$ in the stabilizer $G_p$   such that  
$h\gamma(d)= \delta(d)$. Since the geodesics are small, this implies 
$h\gamma = \delta$ and $\pi (\gamma)= \pi (\delta)$. 
\qed\enddemo 
 
A normal geodesic $\gamma$ in $M$ is called a \idx{\it horizontal 
lift} of the minimal geodesic segment $\pi (\gamma)$ of $M/G$. 
The proof of \therosteritem3 implies that \idx{\it any two horizontal 
lifts of a minimal geodesic segment of $M/G$ differ by the action of 
an isometry} in $G$. 

We define the \idx{\it angle between two minimal geodesic ray 
segments} in $M/G$ from a point $\bar p$ as the minimum of the angles 
between all their horizontal lifts through $p\in \bar p$. The angle is 
independent of the choice of $p$.

\proclaim{\nmb.{3.2}. Proposition }
Let $G$ be a compact connected Lie group  and $\rh: G \to SO(V)$  a
polar orthogonal representation into a Euclidean vector space 
$(V, g= <\quad,\quad>)$. Then any two points of the orbit space $V/G$ are
connected by a unique minimal geodesic segment.
\endproclaim

\demo{Proof} By definition of a polar representation there exists a 
linear section $S$, i.e\., a vector subspace of $V$ which intersects
all orbits orthogonally, see \cit!{8}.
Denote by $W$ the Weyl group of a section $S$,
that is the quotient $W= N_G(S)/Z_G(S)$, where 
$N_G(S)$ is the subgroup of $G$ which preserves $S$ and $Z_G(S)$ is
its normal subgroup which acts trivially on $S$. It is known that
$W$ is a finite group generated by reflections in hyperplanes of $S$
and that the orbit spaces $V/G$ and $S/W$ are isometric, \cit!{8}.
This reduces the statement to the case  when $G$ is a finite group
generated by reflections, i.e.,  the Weyl group of a root system.
The orbit space of such group is the closure of a Weyl chamber
which is a convex polyhedral cone of Euclidean space with the induced 
metric. Now the statement is obvious.
\qed\enddemo

Note that up to now we know only minimal geodesic segments; 
geodesics as we will treat them in \nmb!{4.4} 
below will turn out to be reflected at faces, so there will be many 
different geodesics connecting two points.

\subhead\nmb.{3.3}. Example\endsubhead
Statement \nmb!{3.2} is not true 
if the representation $\rh$ is not polar.
Moreover, there may exist points $x,y$ with arbitrary small distance 
$\ep$ which are connected by several minimal geodesic segments.

Let $G \subset SO(V)$ be a compact connected linear group such that
$V$ is the direct sum of two $G$-invariant subspaces $V_1,V_2$ of
dimension $>1$. Choose $x \in V_1, y \in V_2$ such that the stabilizer
$H=G_x$ does not act transitively on the orbit $G.y$. 
For example, we may assume that $H=G_x=\{\Id\}$. Note that here 
$y\in T_x(G.x)^\bot$, but $y$ is not contained in a slice at $x$. The 
maximal radius of a slice at $x$ is just $|x|$.

Then for any 
$z \in G.y$ the geodesic $\ga(t) =(1-t)x + tz$ is normal, since 
$\ga'(t)=-x+z$ is normal to $X.\ga(t)=(1-t)X.x+tX.z$ for each 
$X\in\frak g\subset \frak s\frak o(V)$. It is a minimal geodesic which 
connects the orbits $G.x$ and $G.y$. 
Hence it defines a minimal geodesic segment $\bar \ga$ in the orbit
space $V/G$. If $z,z'\in Gy$ do not belong to one $H$ orbit, then 
the corresponding geodesics are not $G$-equivalent
and define different minimal geodesic segments of the orbit space 
connecting $\pi(x)$ and $\pi(y)$.

\subhead\nmb.{3.4} \endsubhead
Recall that a subset $N$ of a path metric space $(M,d)$ of   
is called \idx{\it weakly convex} if the induced 
metric on $N$ is also of path metric type. So any two points in $N$ are 
connected in $M$ by some minimal geodesic segment which lies in $N$.

A subset $N$ of a path metric space $(M,d)$    
is called \idx{\it convex} if {\it any} minimal geodesic 
segment in $M$ between two points in $N$ lies in $N$. 

Denote by $(M/G)_{\le(H)}$ the set of orbits with orbit type 
smaller then $(H)$. In particular, $(M/G)_{\le (K)}=  
(M/G)_{\operatorname{reg}}$ if $(K)$ denotes the minimal orbit type.
 
\proclaim{Proposition} $M_{\le(H)}$ is a convex subset of 
$M/G$. In particular, $(M/G)_{\operatorname{reg}}$ is a convex  
open dense submanifold. 
\endproclaim 
  
The proof follows from the following lemma. 

\proclaim{\nmb.{3.5}. Lemma} Let $\bar p \bar q$ be a minimal geodesic  
in $M/G$ of length $d$ and let $pq$ be a horizontal lift. 
Then the stabilizer 
$G_x$ of any interior point of $pq$ is contained in the stabilizers 
$G_p, G_q$ of the end points. 
\endproclaim 

In particular, if $G_p\cap G_q= \{K\}$, the minimal stabilizer group, 
then all interior points of ${\bar p\bar q}$ are regular.

\demo{Proof} Assume for contradiction that there exist an isometry 
$h \in G_x \setminus G_p$. 
Appying $h$ to a minimal geodesic $px$ we obtain a minimal geodesic 
$(hp)x $ of the same length which connects $Gp =\bar p$ and $x$. Then 
the broken geodesic $(hp)xq$ has length $d$ and connects
$\bar p$  and $\bar q$. This is impossible. 
\enddemo 

\head\totoc\nmb0{4}. Vector fields, geodesics, and Jacobi fields on orbit 
spaces \endhead

\subhead\nmb.{4.1}. Smooth vector fields on orbit spaces \endsubhead
Let $M$ be a proper $G$-manifold and let $\Der(C^\infty(M/G))$ denote 
the space of all derivations of the real algebra $C^\infty(M/G)$; 
these are called smooth vector fields on $M/G$. A vector field 
$X\in \Der(C^\infty(M/G))$ is called \idx{\it strata-preserving} if it 
preserves each ideal in $C^\infty(M/G)$ consisting of functions which 
vanish on an orbit stratum. We denote by $\X(M/G)$ the Lie subalgebra of all 
strata-preserving vector fields.

\proclaim{Theorem} \cit!{26} The canonical mapping 
$$
X\mapsto \bar X,\quad
\X(M)^G\to \Der(C^\infty(M/G))
$$ 
has image $\X(M/G)$.
\endproclaim

In \cit!{26} this is stated for compact $G$, but the proof works 
without change also for proper $G$-actions.

Thus we have the following exact sequence:
$$
0 \to \X(M)^G_{\text{ver}} \to \X(M)^G \to \X(M/G) \to 0
$$
where $\X(M)^G_{\text{ver}}$ is the space of all $G$-invariant vector 
fields on $M$ which are tangent to the orbits (vertical). 

Thus a strata preserving smooth vector field $\bar X\in\X(M/G)$ induces a 
derivation on the algebra of smooth functions on each stratum and 
thus a smooth vector field on each stratum which is tangent to this 
stratum. Moreover $\bar X$ induces a local flow on each stratum. By 
using also a lift in $\X(M)^G$,
we get a strata preserving smooth mapping 
$$
\Bbb R\x M/G \supseteq U @>\Fl^{\bar X}>> M/G
$$
which is defined on an open neighborhood $U$ of $0\x M/G$ in 
$\Bbb R\x M/G$. Clearly for $X\in \X(M)^G$ the flows $X$ and 
$\bar X\in \X(M/G)$ are related, $\pi_{M/G}\o\Fl^X_t = 
\Fl^{\bar X}_t\o\pi_{M/G}$, since this is true on each orbit stratum. 

\subhead\nmb.{4.2} \endsubhead
If $M$ is a complete Riemannian $G$-manifold we may also consider the 
vector space $\X(M)^G_{\text{hor}}$ of all smooth \idx{\it horizontal 
$G$-invariant} vector fields, which are normal to each orbit which 
they meet. The space $\X(M)^G_{\text{hor}}$ is a Lie algebra if and 
only if the vector subbundle $\text{Nor}(M)|M_{\text{reg}}$ is 
integrable which is almost equivalent to the fact that $M$ admits a 
section, since then the horizontal bundle over $M_{\text{reg}}$ is 
integrable; there might be topological difficulties, namely, the 
leaves of the horizontal bundle might be not closed. 

\subhead\nmb.{4.3}. Differential equations of second order on $M/G$ 
\endsubhead
Let $\ka:TTM\to TTM$ be the canonical involution. Recall that a 
vector field $\Ga$ on $TM$ is called a
differential equation of second order on $M$ if $\ka\o\Ga=\Ga$ or, 
equivalently, if $T(\pi_M)\o\Ga=\Id_{TM}$. It is called a 
\idx{\it spray} if it 
is quadratic in the sense that 
$\Ga(m^{M}_t.X)=m^{TM}_t.T(m^M_t).\Ga(X)$, where $m^M_t$ is the 
scalar multiplication by $t$ on $TM$. 

Let $M$ be a proper $G$-manifold. Then the extension of the 
$G$-action to $TM$ is also proper. The projection $\bar\Ga$ 
of a $G$-invariant 
second order differential equation or spray $\Ga\in \X(TM)^G$ to 
$\X(TM/G)$ is called a second order differential equation or spray on 
$M/G$. By \nmb!{4.1} we know that $\bar\Ga$ is tangent to all strata 
of $TM/G$, thus integral curves of $\bar\Ga$ make sense and are 
unique as integral curves of vector fields on manifolds. Clearly, the 
projection of an integral curve of $\Ga$ is an integral curve of 
$\bar\Ga$ on $TM/G$.

\subhead\nmb.{4.4}. The geodesic spray on $M/G$ \endsubhead
Fix a complete $G$-invariant Riemannian metric on $M$ and denote by 
$\Ga$ its  geodesic spray. The flow lines of $\Ga$ are the velocity 
fields of geodesics on $M$. The corresponding spray 
$\bar\Ga\in\X(TM/G)$ is called the \idx{\it geodesic spray} on $M/G$. 
Its flow lines are complete and each is contained in one stratum of 
$TM/G$. 

For $x\in M$ let $\Nor(M)_x = (T_x(G.x))^\bot$, the normal space to 
the orbit, which we may split as orthogonal direct sum of the 
subspace $\Nor^{\text{inv}}(M)_x=\Nor(M)_x^{G_x}$ 
which is invariant under the isotropy group $G_x$, and its orthogonal 
complement in $\Nor(M)_x$.
We consider $\Nor(M):=\bigcup_{x\in M}\Nor(M)_x\subset TM$, and 
similarly for $\Nor^{\text{inv}}(M)$.
These are $G$-invariant subsets of $TM$ which can be considered as 
families of 
sub vector spaces with jumping dimensions: Over singular strata the 
dimension may become larger. Note that 
$\Nor(M)_x=\Nor^{\text{inv}}(M)_x$ if and only if $x$ is a regular 
point. 

$\Nor(M)$ is invariant under the 
flow of the spray $\Ga$ since a geodesic which is orthogonal to one 
orbit is orthogonal to any orbit it meets. 
However, $Nor^{inv}(M)$ is not invariant under the flow of
the spray $\Ga$, since a geodesic starting at a regular point 
orthogonally to the orbit may hit later a singular point where its 
tangent vector is still orthogonal to the orbit but no longer 
invariant under the (larger) isotropy group.
Consequently, 
$\Nor(M)/G\subseteq TM/G $ is invariant under the flow of the spray 
$\bar\Ga$. We may consider $\Nor(M)/G$ as a substitute of the
\idx{\it tangent bundle} of the the orbit space $M/G$ since the 
normal slices suffice to describe any tangent vector which moves an 
orbit infinitesimally.

\subhead Definition \endsubhead
A \idx{\it geodesic}  on $M/G$ is a curve of the form
$$
t\mapsto \left(\pi_{M/G}\o \Fl^{\bar \Ga}_t \right)(\bar \xi), 
     \quad  \bar \xi\in \Nor(M)/G
$$
This definition fits well with concept of minimal geodesic arcs as 
treated in section \nmb!{3}; geodesics are prolongations of minimal 
geodesic arcs. Clearly geodesics in $M/G$ are exactly the projections 
onto $M/G$ of normal geodesics in $M$. 

\subhead\nmb.{4.5}. Example \endsubhead
Let $G\to O(V)$ be a polar representation of a compact group, with 
section $\Si\subset V$, a linear subspace which meets every orbit 
orthogonally. Then $V/G=\Si/W(\Si)$ is representated by a chamber $C$ in 
$\Si$. The normal geodesics in $V$ can be chosen to lie in $\Si$, 
thus the geodesics in $V/G=\Si/W(\Si)\cong C$ are straight lines in 
the interior of $C$ which are reflected by the walls:
the incoming angle equals the outgoing angle. 

\subhead\nmb.{4.6}. Ballistic curves \endsubhead
The projection $t\mapsto \pi_{M/G}\o \Fl^{\bar \Ga}_t(\bar \xi)$ onto 
$M/G$ of a flow line of the geodesic spray on $TM/G$ with general initial 
vector $\bar \xi\in TM/G$ (which need \idx{\it not} be in 
$\Nor(M)/G$), is called a \idx{\it ballistic curve}. 
It depends on external data: $\bar\pi_{M/G}:TM/G\to M/G$ is bigger than the 
tangent bundle.

\subhead\nmb.{4.7}. The Jacobian flow \endsubhead
We recall the following result:

\proclaim{Theorem} \cit!{17}
Let $\Ga:TM\to TTM$ be a spray on a manifold $M$. Then 
$\ka_{TM}\o T\Ga:TTM\to TTTM$ is a vector field. Consider a flow line 
$$
Y(t)=\Fl^{\ka_{TM}\o T\Ga}_t(Y(0))
$$
of this field. Then we have:
\roster
\item"" $c:=\pi_M\o\pi_{TM}\o Y$ is a geodesic on $M$.
\item"" $\dot c=\pi_{TM}\o Y$ is the velocity field of $c$.
\item"" $J:=T(\pi_M)\o Y$ is a Jacobi field along $c$.
\item"" $\dot J = \ka_M\o Y$ is the velocity field of $J$.
\item"" $\nabla_{\partial_t}J=K\o \ka_M\o Y$ 
        is the covariant derivative of $J$. 
\item"" The Jacobi equation is given by:
$$\align
0 &= \nabla_{\partial_t}\nabla_{\partial_t}J + R(J,\dot c)\dot c +
     \nabla_{\partial_t}\operatorname{Tor}(J,\dot c)\\
&= K\o TK\o T\Ga\o Y.
\endalign$$
\endroster
This implies that in a canonical chart induced from a chart on $M$ 
the curve $Y(t)$ is given by 
$$
(c(t),c'(t);J(t),J'(t)).
$$
\endproclaim

On a complete Riemannian $G$-manifold $M$ with geodesic spray $\Ga$ 
we may thus consider the the $G$-invariant vector field 
$\ka_{TM}\o T\Ga:TTM\to TTTM$ and the induced smooth derivation
$$
\overline {\ka_{TM}\o T\Ga} \in \X(TTM/G).
$$
We have 
$$
T(\pi_{TM/G})\o \overline {\ka_{TM}\o T\Ga} 
     = \overline{T\pi_{TM}\o \ka_{TM}\o T\Ga}
     = \overline{\pi_{TTM}\o T\Ga}
     = \overline{\Ga\o \pi_{TM}}
     = \bar{\Ga}\o \pi_{TM/G}.
$$
We consider a flow line
$$
\bar Y(t)=\Fl^{\overline{\ka_{TM}\o T\Ga}}_t(\bar Y(0))
$$
of this field which respects the orbit stratification. Then we have:
\roster
\item  $t\mapsto \bar c(t):=\pi_{M/G}\o\pi_{TM/G}\o \bar Y(t)\in M/G$ is a 
       geodesic on $M/G$ if the initial velocity vector 
       $\pi_{TM/G}(\bar Y(0))\in \Nor(M)/G$ is normal. If not then 
       $\bar c(t)$ is a ballistic curve on $M/G$. 
\item  $t\mapsto \dot {\bar c}(t)=\pi_{TM/G}\o \bar Y(t)\in TM/G$ is 
       the velocity field of $c$. It respects the orbit 
       stratification of $TM/G$ since it is a flow line of $\bar \Ga$. 
\item  $t\mapsto \bar J(t):=T(\pi_{M})_{/G}\o \bar Y(t)\in TM/G$ can 
       be called a \idx{\it Jacobi field} along $c$. It does not 
       respect the orbit stratification in general. 
\item  $\dot {\bar J} = \ka_{M/G}\o \bar Y$ is the velocity field of 
       $J$. It respects the orbit stratification. 
\item  $\nabla_{\partial_t}\bar J=\bar K\o \ka_{M/G}\o \bar Y$ 
       is the covariant derivative of $\bar J$. 
\endroster

\head\totoc\nmb0{5}. Example: Hermitian and symmetric matrices \endhead

\subhead\nmb.{5.1}. Simplest Example \endsubhead
Let $S^1$ act on $\Bbb R^2$ by rotations. Then 
$\Bbb R^2/S^1=[0,\infty)$, and 
$\overline{\pi_{\Bbb R^2}}:T\Bbb R^2/S^1\to \Bbb R^2/S^1$ looks as 
follows: The fiber over 0 is $[0,\infty)$ again, and the fiber over 
$t>0$ is $\Bbb R^2$. Normal geodesics on $\Bbb R^2$ are lines 
through 0, so geodesics on $\Bbb R^2/S^1$ are constant speed curves 
coming in from infinity on $[0,\infty)$ which are reflected at 0 and 
go out again at constant speed: $t\mapsto |tv|$.
However, ballistic curves seem to carry a charge and behave as being 
repelled by a field carried by the singular orbit $0$, namely
$$
t\mapsto \sqrt{(x_1+tv_1)^2+(x_2+tv_2)^2},\quad \text{ where }x, v
\text{ are linearly independent in }\Bbb R^2.
$$

\subhead\nmb.{5.2}. Simple Example \endsubhead
Let $SO(2)$ act on the space $S(2)$ of symmetric 
$(2\x 2)$-matrices by conjugation, a polar representation, where the 
diagonal matrices form a section. A chamber is here given by the 
halfspace
$$
S(2)/SO(2)=C:=\left\{\pmatrix \la_1 & 0 \\ 0 & \la_2\endpmatrix: 
     \la_1\ge\la_2\right\}.
$$
Let us describe $TS(2)/SO(2)\to S(2)/SO(2)=C$. Over a point 
$A=\on{diag}(\la,\la)$ in the wall of $C$ we can use the isotropy group 
$SO(2)_A=SO(2)$ to put the tangent vector in normal form, so the 
fiber there is the half space $C$. The fiber over an interior point 
is the whole vector space $S(2)$.

As in \nmb!{4.5} geodesics in $S(2)/SO(2)\cong C$ are 
straight lines which are reflected by the wall $\{\la_1=\la_2\}$. 
One can compute the ballistic curves. 
For 
$$
t\mapsto A+tV = \pmatrix a_1 & 0 \\ 0 & a_2 \endpmatrix + t\;
\pmatrix v_1 & v_3 \\ v_3 & v_2 \endpmatrix \in S(2)
$$
the curve of eigenvalues in $C$ is
$$
t\mapsto \binom{\la_1(t)}{\la_2(t)}= \frac12
\pmatrix a_1+a_2 + t(v_1+v_2) + \sqrt{(a_1-a_2+t(v_1-v_2))^2 + 4t^2v_3^2}\\
         a_1+a_2 + t(v_1+v_2) - \sqrt{(a_1-a_2+t(v_1-v_2))^2 + 4t^2v_3^2}
\endpmatrix.
$$
Here $t\mapsto a_1+a_2 + t(v_1+v_2)$ is the component on the wall 
$\la_1=\la_2$ of $C$ which travels with constant speed, whereas 
$\sqrt{(a_1-a_2+t(v_1-v_2)^2 + 4t^2v_3^2}$ is the distance from the 
wall. So again the ballistic curve is being repelled by a field 
carried by the wall. If it contains one regular orbit and is not a 
geodesic, then it never hits the wall. 

\subhead\nmb.{5.3}. The space of Hermitian matrices \endsubhead
Let $G=SU(n)$ act on the space $H(n)$ of complex Hermitian 
$(n\x n)$-matrices by conjugation, where the inner product is given 
by the (always real) trace $\on{Tr}(AB)$.
This is a polar representation, where the diagonal matrices with real 
entries form a section $\Si$. A chamber is here given by the quadrant 
$C\subset \Si$ consisting of all real diagonal matrices with 
eigenvalues $\la_1\ge\la_2\ge\dots\ge\la_n$.
Geodesics in $H(n)/SU(n)\cong C$ are straight lines which are 
reflected by all walls $\{\la_i=\la_{i_1}=\dots=\la_{i+k}\}$.

A ballistic curve looks as follows: Let $A$ be a diagonal matrix with 
eigenvalues $a_1\ge\dots\ge a_n$, and let $V=(v_{i,j})$ be a 
Hermitian matrix. The ballistic curve is then 
$\la(t)=(\la_1(t)\ge\dots\ge\la_n(t))$, the curve of eigenvalues of 
the Hermitian matrix $A+tV$.

\subhead\nmb.{5.4}. Hamiltonian description \endsubhead
Let us describe ballistic curves as trajectories of a Hamiltonian 
system on a reduced phase space. Let $T^*H(n)=H(n)\x H(n)$ be the 
cotangent bundle where we identified $H(n)$ with its dual by the 
inner product, so the duality is given by 
$\langle \al,A \rangle=\on{Tr}(A\al)$.
Then the canonical 1-form is given by 
$\th(A,\al,A',\al')=\on{Tr}(\al A')$, the symplectic form is
$\om_{(A,\al)}((A',\al'),(A'',\al''))=\on{Tr}(A'\al''-A''\al')$, and 
the Hamiltonian function for the straight lines $(A+t\al,\al)$ on 
$H(n)$ is $h(A,\al)=\frac12\on{Tr}(\al^2)$. The action 
$SU(n)\ni g\mapsto (A\mapsto gAg\i)$ lifts to the action 
$SU(n)\ni g\mapsto ((A,\al)\mapsto (gAg\i,g\al g\i))$ on $T^*H(n)$ 
with fundamental vector fields $\ze_X(A,\al)=(A,\al,[X,A],[X,\al])$ 
for $X\in \frak s\frak u(n)$, and with generating functions 
$f_X(A,\al)=\th(\ze_X(A,\al))=\on{Tr}(\al[X,A])=\on{Tr}([A,\al]X)$. 
Thus the momentum mapping  $J:T^*H(n)\to \frak s\frak u(n)^*$ is given by 
$\langle X,J(A,\al) \rangle=f_X(A,\al)= \on{Tr}([A,\al]X)$. If we 
identify $\frak s\frak u(n)$ with its dual via the inner product 
$\on{Tr}(XY)$, the momentum mapping is $J(A,\al)=[A,\al]$.
Along the line $t\mapsto A+t\al$ the momentum mapping is constant: 
$J(A+t\al,\al)=[A,\al]=Y\in \frak s\frak u(n)$.
Note that for $X\in \frak s\frak u(n)$ the evaluation on $X$ of 
$J(A+t\al,\al)\in \frak s\frak u(n)^*$ equals the inner product:  
$$
\langle  X,J(A+t\al,\al)\rangle = 
\on{Tr}(\tfrac d{dt}(A+t\al),\ze_X(A+t\al)),
$$
which is obviously constant in $t$; compare with the general 
result of Riemannian transformation groups, e.g.\ \cit!{16},~8.1. 

According to principles of symplectic reduction 
\cit!{1},~4.3.5, \cit!{27}, \cit!{12}, \cit!{22}, we have to 
consider for a regular value $Y$ (and later for an arbitrary value) 
of the momentum mapping $J$ the 
submanifold $J\i(Y)\subset T^*H(n)$. The null distribution of 
$\om|J\i(Y)$ is integrable (with jumping dimensions) and its leaves 
(according to the Stefan-Sussmann theory of integrable distributions)
are exactly the orbits in $J\i(Y)$ of the isotropy group $SU(n)_Y$ for the 
coadjoint action. 
So we have to consider the orbit space $J\i(Y)/SU(n)_Y$.
If $Y$ is not a regular value of $J$, the inverse image $J\i(Y)$ is a 
subset which is described by polynomial equations since $J$ is 
polynomial (in fact quadratic), so $J\i(Y)$ is stratified into 
submanifolds; symplectic reduction 
works also for this case, see \cit!{27}.

\subhead\nmb.{5.5}. The case of momentum $Y=0$ gives again geodesics 
\endsubhead If $Y=0$ then $SU(n)_Y=SU(n)$ and 
$J\i(0)=\{(A,\al):[A,\al]=0\}$, so $A$ and $\al$ commute. If $A$ is 
regular (i.e. all eigenvalues are distinct), using a uniquely 
determined transformation $g\in SU(n)$ we move the point $A$ into the 
open chamber $C^o\subset H(n)$, so $A=\on{diag}(a_1> a_2>\dots> a_n)$ 
and since $\al$ commutes with $A$ so it is also in diagonal form. The 
symplectic form $\om$ restricts to the canonical symplectic form on 
$C^o\x \Si = C^o\x \Si^* = T^*(C^o)$. Thus symplectic reduction gives
$(J\i(0)\cap (T^*H(n))_{\text{reg}})/SU(n)= T^*(C^o)\subset T^*H(n)$. 
By \cit!{27} we also use symplectic reduction for non-regular $A$ and 
we get (see in particular \cit!{12},~3.4) $J\i(0)/SU(n)=T^*C$, the 
stratified cotangent cone bundle of the chamber $C$ considered 
asstratified space. Namely, if one root $\ep_i(A)=a_i-a_{i+1}$ 
vanishes on the diagonal matrix $A$ then the isotropy group $SU(n)_A$ 
contains a subgroup $SU(2)$ corresponding to these coordinates. Any 
matrix $\al$ with $[A,\al]=0$ contains an arbitrary hermitian 
submatrix corresponding to the coordinates $i$ and $i+1$, which may 
be brougth into diagonal form with the help of this $SU(2)$ so that 
$\ep_i(\al)=\al_i-\al_{i+1}\ge 0$. Thus the tangent vector $\al$ with 
foot point in a wall is either tangent to the wall (if 
$\al_i=\al_{i+1}$) or points into the interior of the chamber $C$. 
The Hamiltonian $h$ restricts to 
$C^o\x\Si\ni (A,\al)\mapsto \frac12\sum_i\al_i^2$, so the 
trajectories of the Hamiltonian system here are again straight lines 
which are reflected at the walls. 

\subhead\nmb.{5.6}. The case of general momentum $Y$ \endsubhead
If $Y\ne 0\in \frak s\frak u(n)$ and if $SU(n)_Y$ is the isotropy group 
of $Y$ for the adjoint representation, then it is well known (see 
references in \nmb!{5.4}) that we 
may pass from $Y$ to the coadjoint orbit $\Cal O(Y)=\on{Ad}^*(SU(n))(Y)$ 
and get 
$$
J\i(Y)/SU(n)_Y = J\i(\Cal O(Y))/SU(n) = (J\i(Y)\x \Cal O(-Y))/SU(n),
\tag{\nmb:{1}}$$
where all (stratified) diffeomorphisms are symplectic ones.

\subhead\nmb.{5.7}. The Calogero Moser system  \endsubhead
As the simplest case we assume that $Y'\in \frak s\frak u(n)$ is 
not zero but has maximal isotropy group, see \cit!{10}: 
So we assume that $Y'$ has 
complex rank 1 plus an imaginary multiple of the identity, 
$Y'=\sqrt{-1}(c\Bbb I_n + v\otimes v^*)$ for $0\ne v=(v^i)$ a column 
vector in $\Bbb C^n$. The coadjoint orbit is then 
$\Cal O(Y')=\{\sqrt{-1}(c\Bbb I_n + w\otimes w^*):w\in \Bbb C^n, 
|w|=|v|\}$, isomorphic to $S^{2n-1}/S^1=\Bbb CP^n$, of real 
dimension $2n-2$.  
Consider $(A',\al')$ with  $J(A',\al')=Y'$, choose $g\in SU(n)$ 
such that $A=gA'g\i=\on{diag}(a_1\ge a_2\ge \dots\ge a_n)$, and let 
$\al=g\al'g\i$. 
Then the entry of the commutator is $[A,\al]_{ij}=\al_{ij}(a_i-a_j)$.
So $[A,\al]=gY'g\i =:Y= \sqrt{-1}(c\Bbb I_n + gv\otimes (gv)^*) 
= \sqrt{-1}(c\Bbb I_n + w\otimes w^*)$
has zero diagonal entries, thus 
$0<w^i\bar w^i=-c$ and $w^i=\exp(\sqrt{-1}\th_i)\sqrt{-c}$ for some 
$\th_i$
But then all off-diagonal entries $Y_{ij}=\sqrt{-1}w^i\bar w^j = 
-\sqrt{-1}\,c\,\exp(\sqrt{-1}(\th_i-\th_j))\ne 0$, and $A$ has to be 
regular. 
We may use the remaining gauge freedom in the isotropy group 
$SU(n)_A=S(U(1)^n)$ to put $w^i=\exp(\sqrt{-1}\th)\sqrt{-c}$ where 
$\th=\sum\th_i$. Then $Y_{ij}=-c\sqrt{-1}$ for $i\ne j$. 

So the reduced space $(T^*H(n))_{Y}$ is diffeomorphic to the 
submanifold of $T^*H(n)$ consisting of all 
$(A,\al)\in H(n)\x H(n)$ where $A=\on{diag}(a_1>a_2>\dots>a_n)$, 
and where $\al$ has arbitrary diagonal entries $\al_i:=\al_{ii}$ and 
off-diagonal entries $\al_{ij}=Y_{ij}/(a_i-a_j)=-c\sqrt{-1}/(a_i-a_j)$.
We can thus use $a_1,\dots,a_n,\al_1,\dots,\al_n$ as coordinates. 
The invariant symplectic form pulls back to 
$\om_{(A,\al)}((A'\al'),(A'',\al''))=\on{Tr}(A'\al''-A''\al')
     =\sum(a_i'\al_i''-a_i''\al_i')$.
The invariant Hamiltonian $h$ restricts to the Hamiltonian
$$
h(A,\al) = \tfrac12\on{Tr}(\al^2)
 = \frac12\sum_i \al_i^2 +\frac12\sum_{i\ne j}\frac{c^2}{(a_i-a_j)^2}.
$$
This is the famous Hamiltonian function of the Calogero-Moser 
completely integrable system, see \cit!{18}, \cit!{19}, \cit!{10}, 
and \cit!{22},~3.1 and 3.3. The corresponding Hamiltonian vector 
field and the differential equation for the ballistic curve are then 
$$\gather
H_h=\sum_i \al_i\frac{\partial}{\partial a_i} +
     2\sum_i\sum_{j\ne i}\frac{c^2}{(a_i-a_j)^3}
     \frac{\partial}{\partial \al_i},\\
\ddot a_i = 2\sum_{j\ne i}\frac{c^2}{(a_i-a_j)^3},\\
\endgather$$
Note that the ballistic curve avoids the walls of the Weyl chamber $C$.

\subhead\nmb.{5.8}. Degenerate cases of non-zero momenta of minimal 
rank \endsubhead
Let us discuss now the case of non-regular diagonal $A$. 
Namely, if one root, say $\ep_{12}(A)=a_1-a_2$ vanishes on the 
diagonal matrix $A$ then the isotropy group $SU(n)_A$ contains a 
subgroup $SU(2)$ corresponding to these coordinates. 
Consider $\al$ with $[A,\al]=Y$; 
then $0=\al_{12}(a_1-a_2)=Y_{12}$. 
Thus $\al$ 
contains an arbitrary hermitian submatrix corresponding to the 
first two coordinates, which may be brougth into diagonal form 
with the help of this $SU(2)\subset SU(n)_A$ so that 
$\ep_{12}(\al)=\al_1-\al_2\ge 0$. Thus the tangent vector $\al$ with 
foot point $A$ in a wall is either tangent to the wall (if 
$\al_1=\al_2$) or points into the interior of the chamber $C$ 
(if $\al_1>\al_2$). 
Note that then $Y_{11}=Y_{22}=Y_{12}=0$.

Let us now assume that the momentum $Y$ is of the form 
$Y=\sqrt{-1}(c\Bbb I_{n-2} + v\otimes v^*)$ for some vector 
$0\ne v\in \Bbb C^{n-2}$. We can repeat the analysis of \nmb!{5.7}
in the subspace $\Bbb C^{n-2}$, and get for the Hamitonian
$$\gather
h(A,\al) = \tfrac12\on{Tr}(\al^2) = \frac12\sum_{i=1}^n \al_i^2 
      +\frac12\sum_{3\le i\ne j}\frac{c^2}{(a_i-a_j)^2},\\
H_h=\sum_{i=1}^n \al_i\frac{\partial}{\partial a_i} +
     2\sum_{3\le i\ne j}\frac{c^2}{(a_i-a_j)^3}
     \frac{\partial}{\partial \al_i},\\
\ddot a_1=\ddot a_2=0,\quad 
     \ddot a_i = 2\sum_{3\le j\ne i}\frac{c^2}{(a_i-a_j)^3}\text{ for 
     }i>2. 
\endgather$$
So the ballistic curves are just the trajectories of the 
Calogero-Moser integrable system inside the wall $\{a_1-a_2=0\}$ 
complemented by a geodesic in the coordinates orthogonal to this wall. 
Of course we may add other vanishing roots. 

\subhead\nmb.{5.9}. The case of general momentum $Y$ and regular $A$ 
\endsubhead
Starting again with some regular $A'$ consider 
$(A',\al')$ with  $J(A',\al')=Y'$, choose $g\in SU(n)$ 
such that $A=gA'g\i=\on{diag}(a_1>a_2>\dots>a_n)$, and let 
$\al=g\al'g\i$ and $Y=gY'g\i=[A,\al]$.
Then the entry of the commutator is 
$Y_{ij}=[A,\al]_{ij}=\al_{ij}(a_i-a_j)$ thus $Y_{ii}=0$. We may pass 
to the coordinates $a_i$ and $\al_i:=\al_{ii}$ for $1\le i\le n$ on 
the one hand, corresponding to $J\i(Y)$ in \nmb!{5.6.1}, and 
$Y_{ij}$ for $i\ne j$ on the other hand, corresponding to $\Cal O(-Y)$ in 
\nmb!{5.6.1}, with the linear relation 
$Y_{ji}=-\overline{Y_{ij}}$ and with $n-1$ non-zero entries 
$Y_{ij}>0$ with $i>j$
(chosen in lexicographic order) by 
applying the remaining isotropy group 
$SU(n)_A=S(U(1)^n)=\{\on{diag}(e^{\sqrt{-1}\th_1},\dots,e^{\sqrt{-1}\th_n}):
     \sum\th_i\in 2\pi\Bbb Z\}$.
We may use this canonical form as section 
$$(J\i(Y)\x \Cal O(-Y))/SU(n)\to J\i(Y)\x \Cal O(-Y)
     \subset TH(n)\x \frak s\frak u(n)$$ 
to pull back the symplectic or Poisson structures and the Hamiltonian 
function
$$\align
h(A,\al) &= \tfrac12\on{Tr}(\al^2) 
= \frac12\sum_i \al_i^2 
      -\frac12\sum_{i\ne j}\frac{Y_{ij}Y_{ji}}{(a_i-a_j)^2},\\
dh &= \sum_i \al_i\,d\al_i 
      +\sum_{i\ne j}\frac{Y_{ij}Y_{ji}}{(a_i-a_j)^3}(da_i-da_j)
      -\frac12\sum_{i\ne j}\frac{dY_{ij}.Y_{ji}+Y_{ij}.dY_{ji}}{(a_i-a_j)^2},\\
&= \sum_i \al_i\,d\al_i 
      +2\sum_{i\ne j}\frac{Y_{ij}Y_{ji}}{(a_i-a_j)^3}da_i
      -\sum_{i\ne j}\frac{Y_{ji}}{(a_i-a_j)^2}dY_{ij}.\tag{\nmb:{1}}\\
\endalign$$
The invariant symplectic form on $TH(n)$ pulls back to 
$\om_{(A,\al)}((A'\al'),(A'',\al''))=\on{Tr}(A'\al''-A''\al')
     =\sum(a_i'\al_i''-a_i''\al_i')$ 
thus to $\sum_ida_i\wedge d\al_i$. 
The Poisson structure on $\frak s\frak u(n)$ is given by 
$$\allowdisplaybreaks\align
\La_Y(U,V)&=\on{Tr}(Y[U,V])
     =\sum_{m,n,p}(Y_{mn}U_{np}V_{pm}-Y_{mn}V_{np}U_{pm})\\
\La_Y &= \sum_{i\ne j, k\ne l}
     \La_Y(dY_{ij},dY_{kl})\partial_{Y_{ij}}\otimes\partial_{Y_{kl}}\\
&= \sum_{i\ne j, k\ne l}\sum_{m,n}
     (Y_{mn}\de_{ni}\de_{jk}\de_{lm}-Y_{mn}\de_{nk}\de_{li}\de_{jm})
     \partial_{Y_{ij}}\otimes\partial_{Y_{kl}}\\
&= \sum_{i\ne j, k\ne l}
     (Y_{li}\de_{jk}-Y_{jk}\de_{li})
     \partial_{Y_{ij}}\otimes\partial_{Y_{kl}}\\
\endalign$$
Since this Poisson 2-vector field is tangent to the orbit 
$\Cal O(-Y)$ and is $SU(n)$-invariant, we can push it down to the 
orbit space. There it maps $dY_{ij}$ to (remember that $Y_{ii}=0$)
$$
\La_{-Y}(dY_{ij})
     = -\sum_{k\ne l}(Y_{li}\de_{jk}-Y_{jk}\de_{li})\partial_{Y_{kl}}
     = -\sum_{k}(Y_{ki}\partial_{Y_{jk}}-Y_{jk}\partial_{Y_{ki}}).
$$
So by \thetag{\nmb|{1}} the Hamiltonian vector field is
$$\align
H_h&= \sum_i \al_i\,\partial_{a_i} 
      -2\sum_{i\ne j}\frac{Y_{ij}Y_{ji}}{(a_i-a_j)^3}\,\partial_{\al_i}
      +\sum_{i\ne j}\frac{Y_{ji}}{(a_i-a_j)^2}
      \sum_{k}(Y_{ki}\,\partial_{Y_{jk}}-Y_{jk}\,\partial_{Y_{ki}})\\
&= \sum_i \al_i\,\partial_{a_i} 
      -2\sum_{i\ne j}\frac{Y_{ij}Y_{ji}}{(a_i-a_j)^3}\,\partial_{\al_i}
      -\sum_{i,j,k}\left(\frac{Y_{ji}Y_{jk}}{(a_i-a_j)^2}
      -\frac{Y_{ij}Y_{kj}}{(a_j-a_k)^2}\right)\partial_{Y_{ki}}\\
\endalign$$
The differential equation thus becomes (remember that $Y_{jj}=0$):
$$\align
\dot a_i &= \al_i \\  
\dot\al_i &= -2\sum_{j}\frac{Y_{ij}Y_{ji}}{(a_i-a_j)^3} 
     = 2\sum_{j}\frac{|Y_{ij}|^2}{(a_i-a_j)^3}\\ 
\dot Y_{ki} &= -\sum_{j}\left(\frac{Y_{ji}Y_{jk}}{(a_i-a_j)^2}
      -\frac{Y_{ij}Y_{kj}}{(a_j-a_k)^2}\right).
\endalign$$
Consider the Matrix $Z$ with $Z_{ii}=0$ and 
$Z_{ij}=Y_{ij}/(a_i-a_j)^2$. Then the differential equations become:
$$
\ddot a_i = 2\sum_{j}\frac{|Y_{ij}|^2}{(a_i-a_j)^3},\qquad
\dot Y = [Y^*,Z].
$$
This is the Calogero-Moser integrable system with spin, see \cit!{3} 
and \cit!{4}. 

\subhead\nmb.{5.10}. The case of general momentum $Y$ and singular $A$ 
\endsubhead
Let us consider the situation of \nmb!{5.9}, when $A$ is not regular. 
Let us assume again that one root, say $\ep_{12}(A)=a_1-a_2$ vanishes 
on the  
diagonal matrix $A$. Consider $\al$ with $[A,\al]=Y$. From 
$Y_{ij}=[A,\al]_{ij}=\al_{ij}(a_i-a_j)$ we conclude that $Y_{ii}=0$ 
for all $i$ and also $Y_{12}=0$. The isotropy group $SU(n)_A$ 
contains a subgroup $SU(2)$ corresponding to the first two 
coordinates and we may use this to move $\al$ into the form that 
$\al_{12}=0$ and $\ep_{12}(\al)\ge0$. Thus the tangent vector $\al$ 
with foot point $A$ in the wall $\{\ep_{12}=0\}$ is either tangent to 
the wall when $\al_1=\al_2$ or points into the interior of the 
chamber $C$ when $\al_1>\al_2$. We can then use the same analysis as 
in \nmb!{5.9} where we use now that $Y_{12}=0$. 

In the general case, when some roots vanish, we get for the 
Hamiltonian function, vector field, and differential equation:
$$\align
h(A,\al) &= \tfrac12\on{Tr}(\al^2) = \frac12\sum_i \al_i^2 
      +\frac12\sum_{\{(i,j):a_i\ne a_j\}}\frac{|Y_{ij}|^2}{(a_i-a_j)^2},\\
H_h&=\sum_i \al_i\partial_{a_i} +
     2\sum_{a_j\ne a_i}\frac{|Y_{ij}|^2}{(a_i-a_j)^3}
     \,\partial_{\al_i} + \\
&\quad -\sum_{a_i\ne a_j,k}\frac{Y_{ji}Y_{jk}}{(a_i-a_j)^2}\partial_{Y_{ki}}
     +\sum_{i,a_j\ne a_k}\frac{Y_{ij}Y_{kj}}{(a_j-a_k)^2}\partial_{Y_{ki}}\\
\ddot a_i &= 2\sum_{a_j\ne a_i}\frac{|Y_{ij}|^2}{(a_i-a_j)^3},\qquad
\dot Y = [Y^*,Z]
\endalign$$
where we use the same notation as above. It would be very interesting 
to investigate the reflection behavior of this ballistic curve at the 
walls. 

\subhead\nmb.{5.11}. Example: symmetric matrices \endsubhead
We finally treat the action of $SO(n)=SO(n,\Bbb R)$ on the space 
$S(n)$ of symmetric matrices by conjugation. 
Following the method of 
\nmb!{5.9} and \nmb!{5.10} we get the following result. 
Let $t\mapsto A'+t\al'$ be a straight line in $S(n)$. Then the 
ordered set of eigenvalues $a_1(t),\dots,a_n(t)$ of $A'+t\al'$
is part of the integral curve of the following vector field:
$$\align
H_h&=\sum_i \al_i\partial_{a_i} +
     2\sum_{a_j\ne a_i}\frac{Y_{ij}^2}{(a_i-a_j)^3}
     \,\partial_{\al_i} + \\
&\quad +\sum_{a_i\ne a_j,k}\frac{Y_{ij}Y_{jk}}{(a_i-a_j)^2}\partial_{Y_{ki}}
     -\sum_{i,a_j\ne a_k}\frac{Y_{ij}Y_{jk}}{(a_j-a_k)^2}\partial_{Y_{ki}}\\
\ddot a_i &= 2\sum_{a_j\ne a_i}\frac{Y_{ij}^2}{(a_i-a_j)^3},\qquad
\dot Y = [Y,Z],\qquad\text{ where }Z_{ij}=-\frac{Y_{ij}}{(a_i-a_j)^2},
\endalign$$
where we also note that $Y_{ij}=Z_{ij}=0$ whenever $a_i=a_j$.

\subhead\nmb.{5.12}. Remark \endsubhead
Along the same line one can investigate the action of the 
quaternionic unitary group $Sp(n)$ on the space of all quaternionic 
hermitian matrices. 
Since the results are more complicated to write down and since they 
are a special case of section~\nmb!{6} we do not dwell on them here.

\head\totoc\nmb0{6}. Ballistic curves on polar representations \endhead

\subhead\nmb.{6.1}. The setting \endsubhead
Let $\rh:G\to O(V,\langle\quad,\quad \rangle)$ be a polar 
representation of a compact connected semisimple group, see 
\nmb!{4.5}, with section $\Si\subset V$, a linear subspace which 
meets every orbit orthogonally. Then $V/G=\Si/W(\Si)$ is 
representated by a chamber $C$ in $\Si$. The normal geodesics in $V$ 
can be chosen to lie in $\Si$, thus the geodesics in 
$V/G=\Si/W(\Si)\cong C$ are straight lines in the interior of $C$ 
which are reflected by the walls.

By Dadok, \cit!{8},~proposition~6, which follows from his 
classification, for any polar representation 
there exists an isotropy representation of a symmetric space with the 
same orbits, and it suffices to investigate those latter ones. 
Thus we can assume that 
$\frak l = \g \oplus V$ is a reductive decomposition of a compact 
semisimple  
Lie algebra $\frak l$, where $\g$ is the compact Lie algebra of $G$, 
and where $\langle\quad,\quad \rangle$ is an invariant positive 
definite inner product on $\frak l$, the negative of the Killing 
form, and where 
$\frak l=\g\oplus V$ is an orthogonal decomposition.
Moreover the infinitesimal action $\rh'$ of $\g$ on 
$V$ is by the adjoint action, $\rh'(X)A=[X,A]$, so $[\g,V]\subseteq V$, and 
$[V,V]\subseteq \g$. 
As section $\Si$ we may use any maximal abelian 
subspace in $V$. Moreover we shall use the following lemma.

\proclaim{\nmb.{6.2}. Lemma} In the situation above we have:
\roster
\item "(\nmb:{1})" 
       For $A\in V$ let $\g'$ be the orthogonal complement of the 
       centralizer $Z_{\g}(A)$ in $\g$ and let $V'$ be the  
       orthogonal complement to the centralizer $Z_V(A)$ in $V$ such 
       that $\frak l= Z_{\g}(A)\oplus \g'\oplus Z_V(A)\oplus V'$. 

       Then $\ad_A$ induces linear isomorphisms 
       $\ad_A:V'\to [A,V']=\g'$ and 
       $\ad_A:\g'\to [A,\g']=V'$.
\item "(\nmb:{2})" 
       An element $A\in V$ is regular for the $G$-action 
       on $V$ if and only if the centralizer $Z_V(A)$ in $V$ is a 
       maximal commutative subalgebra of $V$. In this case $Z_V(A)$ 
       is the unique section in $V$ containing $A$.
\item "(\nmb:{3})" 
       An element $A\in V$ is regular for the $G$-action if and only 
       if there exists an element $X\in\g$ such that $X+A$ is regular 
       in $\frak l$, so that $Z_{\frak l}(X+A)$ is a Cartan 
       subalgebra of $\frak l$.
\item "(\nmb:{4})" 
       A linear subspace $\Si\subset V$ is a section if and only if 
       there exists a Cartan subalgebra $\frak h$ of $\g$ such that 
       $\frak h\oplus \Si$ is a Cartan subalgebra of $\frak l$.  

       In this case, let $R\subset L(\Si,\Bbb R)$ be the system of 
       restricted roots so that we have the orthogonal root space 
       decomposition
$$
\frak l = \frak h \oplus \Si \oplus \bigoplus_{\la\in R}\frak l_\la
$$
       where each $\frak l_\la$ has an orthogonal basis 
       $E_\la^i,B_\la^i$, where $i=1,\dots,k_\la$, and where 
       $E_\la^i\in\g$ and $B_\la^i\in V$ are unit vectors, such that 
       $[A,E_\la^i]=\la(A)B_\la^i$ and $[A,B_\la^i]=\la(A)E_\la^i$ 
       for all $A\in \Si$.  
\item "(\nmb:{5})" 
       Let $\Si\subset V$ be a section for the $G$-action on $V$ and 
       let $A\in \Si$. Then for any $\al\in V$ with $[A,\al]=0$ there 
       exists then there exists some $g\in G_A=Z_G(A)$ with 
       $g.\al\in \Si$. 
\endroster
\endproclaim

\demo{Proof}
\therosteritem{\nmb|{1}}
The result follows from the Fitting decomposition 
$\frak l = Z_{\frak l}(A) + \ad_A(\frak l)$ of the skew-symmetric 
endomorphism $\ad_A:\g\to\g$ which we may write as 
$$
\frak l = (Z_{\g}(A)\oplus \ad_A(V'))\oplus (Z_V(A)\oplus \ad_A(\g'))
$$
since $\ad_A$ interchanges $\g$ and $V$. 

\therosteritem{\nmb|{2}} It is known (see e.g.\ \cit!{8}) that the 
$G$-regular elements in a polar $G$-module $V$ are exactly those 
$A\in V$ such that $Z_{\g}(A)$ is of minimal dimension. Moreover, any 
$A\in V$ is contained in some section $\Si$, and any section $\Si$ is 
a maximal commutative subspace of $V\subset \frak l=\g\oplus V$. 
Since a section $\Si$ generates a compact torus in the Lie group $L$, 
there exists an element $A\in \Si$ such that $Z_V(A)=\Si$. These 
elements are characterized as those such that $Z_V(A)$ is of minimal 
dimension, since for any $A\in V$ the centralizer $Z_A(V)$ contains a 
section $\Si$. By \therosteritem{\nmb|{1}} the element $A\in V$ is 
$G$-regular if and only if $\dim(Z_V(A))$ is minimal, i.e.\ 
$Z_V(A)=\Si$ is a section. 

\therosteritem{\nmb|{3}} follows from \therosteritem{\nmb|{2}}.

\therosteritem{\nmb|{4}}  The first assertion follows from 
\therosteritem{\nmb|{3}}, and the rest is well known in the theory of 
symmetric spaces.

\therosteritem{\nmb|{5}} Denote by $\Si'$ a maximal commutative 
subspace of $V$ containing $A$ and $\al$. Then $\Si'$ is a section. 
Now the isotropy group $G_A=Z_G(A)$ acts transitively on the set of 
all sections of $V$ which contain $A$. So in particular there exists 
$g\in G_A$ such that $g.\Si'=\Si$, and hence $g.\al\in\Si$. 
\qed\enddemo

\subhead\nmb.{6.3}. Hamiltonian description and symplectic reduction
\endsubhead
Generalizing \nmb!{5.4}, 
under the assumption of \nmb!{6.1}, we consider 
$TV=V\x V\cong V\x V^*=T^*V$ where we identify $V$ with $V^*$ via the 
inner product. The canonical 1-form is then given by 
$\th(A,\al,A',\al')=\langle\al, A'\rangle$, the symplectic form is 
$\om_{(A,\al)}((A',\al'),(A'',\al''))=\langle A',\al'' 
\rangle-\langle A''\al' \rangle$,
the Hamiltonian function for the straight lines $(A+t\al,\al)$ is 
$h(A,\al)=\frac12\langle \al,\al \rangle$. The action of $G$ on $V$ 
lifts to the diagonal action on $TV=V\x V$ with fundamental vector 
field $\ze_X(A,\al)=(A,\al,[X,A],[X,\al])$ for $X\in \g$ and with 
generating function $f_X(A,\al)=\th(\ze_X(A,\al))=\langle \al,[X,A] 
\rangle=\langle X,[A,\al] \rangle$, where we used the invariance of 
the inner product on $\frak l$. Thus the momentum mapping 
$J:TV\to \g^*\cong \g$ for this action is given by 
$\langle  X,J(A,\al)\rangle = f_X(A,\al)$, so $J(A,\al)=[A,\al]$. 
Along each line the momentum is constant, $J(A+t\al,\al)=[A,\al]$.

According to principles of symplectic reduction 
\cit!{1},~4.3.5, \cit!{27}, \cit!{12}, \cit!{22}, we have to 
consider a regular value $Y\in\g$ (and later for an arbitrary value) 
of the momentum mapping $J$ and the 
submanifold $J\i(Y)\subset TV$. The null distribution of 
$\om|J\i(Y)$ is integrable (with jumping dimensions) and its leaves 
(according to the Stefan-Sussmann theory of integrable distributions)
are exactly the orbits in $J\i(Y)$ of the isotropy group $G_Y$ for the 
coadjoint action. 
So we have to consider the orbit space $J\i(Y)/G_Y$.
If $Y$ is not a regular value of $J$, the inverse image $J\i(Y)$ is a 
subset which is described by polynomial equations since $J$ is 
polynomial (in fact quadratic), so $J\i(Y)$ is stratified into 
submanifolds; symplectic reduction 
works also for this case, see \cit!{27}. 

\subhead\nmb.{6.4}. The case of momentum $Y=0$ gives again geodesics 
\endsubhead 
If $Y=0$ then $G_Y=G$ and $J\i(0)=\{(A,\al):[A,\al]=0\}$, so 
$A$ and $\al$ commute. We may use $g\in G$ to move $(A,\al)$ such 
that 
$A\in C\subset\Si\subset V$. If $A\in\Si$ is regular for the 
$G$-action then by lemma \nmb!{6.2.2} it is also regular in the sense that 
$\Si=Z_{V}(A)$, thus $\al\in \Si$ also, and the reduced 
phase space is $TC^o=C^o\x \Si$. 
The Hamiltonian is still $h(A,\al)=\frac12\langle  \al,\al\rangle$,
the line $A+t\al$ is in $\Si$, and thus projects to a geodesic in $C$ 
which is reflected at walls. 

If $A$ is not regular and $[A,\al']=0$ then there exists $g$ in the 
isotropy group $G_A$ such that $\al=g.\al'\in C$, by \nmb!{6.2.5}. 
If $A$ is in a wall $\{\la=0\}$, and if $\la(\al)=0$, then the 
geodesic is in the same wall. If $\la(\al)>0$ then the geodesic is 
reflected at this wall. 

\subhead\nmb.{6.5}. Symplectic reduction at regular points\endsubhead
If $Y\ne 0\in \g$ and if $G_Y$ is the isotropy group 
of $Y$ for the adjoint representation, then it is well known that we 
may pass from $Y$ to the adjoint orbit $\Cal O(Y)=\on{Ad}^*(G)(Y)$ 
and get 
$$
J\i(Y)/G_Y = J\i(\Cal O(Y))/G = (J\i(Y)\x \Cal O(-Y))/G,
\tag{\nmb:{1}}$$
where all (stratified) diffeomorphisms are symplectic ones.

We start again with some regular $A\in V$ which we may move into 
$C^o\subset\Si\subset V$ by using some suitable $g\in G$.
Given $\al\in V$ we consider $[A,\al]=Y\in\g$. According to the 
restricted root space decomposition \nmb!{6.2}
we can decompose $\al=\al_\Si+\sum_{\la\in R, i}\al_\la^iB_\la^i$ where 
$\al_\Si\in\Si$. Similarly we 
decompose $Y=Y_{\frak h}+\sum_{\la\in R,i}Y_\la^i.E_\la^i$. Then  
$Y=[A,\al]$ implies $Y_{\frak h}=0\in \frak h$ and 
$Y_\la^i=\la(A)\al_\la^i$. 
Since $A$ is regular, $A\in C^o$ and thus 
$\la(A)\ne 0$ for all $\la\in R$.
Let us use also an orthonormal basis $B_0^i$ of $\Si$ to expand 
$\al_\Si=\sum_i\al_0^iB_0^i$ and $A=\sum_iA_0^iB_0^i$.

We can thus use as coordinates
$$\align
&A_0^i,\al_0^i\in C^o\x \Si = TC^o\\
&Y_\la^i=\la(A)\al_\la^i\in \on{Ad}(G).Y 
\endalign$$
The Hamiltonian function in these splitting is given by
$$\align
h(A,\al) &= \frac12\langle \al,\al\rangle 
     = \frac12\sum_i(\al_0^i)^2 + 
     \frac12\sum_{\la\in R,i}\frac{(Y_\la^i)^2}{\la(A)^2} \\
dh &= \sum_i\al_0^i\,d\al_0^i -
\sum_{\la\in R,i}\frac{(Y_\la^i)^2}{\la(A)^3}\sum_k\la(B_0^k)\,dA_0^k
+ \sum_{\la\in R,i}\frac{Y_\la^i}{\la(A)^2}dY_\la^i \\
\endalign$$
The Poisson 2-field on $\g$ (which is tangent to each adjoint orbit) 
is given by 
$$\align
\La_Y(dU,dV) &= \langle Y,[U,V] \rangle,\quad Y,U,V\in\g\\
\La_Y&=\sum_{\la,\mu,i,j}\La_Y(dY_\la^i,dY_\mu^j)
     \partial_{Y_\la^i}\otimes\partial_{Y_\mu^j} \\
&=\sum_{\la,\mu,i,j}\langle Y,[E_\la^i,E_\mu^j]\rangle
     \partial_{Y_\al^i}\otimes\partial_{Y_\mu^j} \\
&=\sum_{\la,\mu,i,j}\langle Y,\sum_k N^{ijk}_{\la\mu} E_{\la+\mu}^k\rangle
     \partial_{Y_\al^i}\otimes\partial_{Y_\mu^j} \\
\La_{-Y}(dY_\la^i) 
     &=\sum_{\mu,j}\langle -Y,\sum_k N^{ijk}_{\la\mu} E_{\la+\mu}^k\rangle
     \partial_{Y_\mu^j} \\
\endalign$$
where we used the convention 
$[E^i_\la,E^j_\mu]=\sum_k N^{ijk}_{\la\mu}E_{\la+\mu}^k$.
We get the following Hamiltonian vector field
$$\align
H_h &= \sum_i\al_0^i\,\partial_{A_0^i} +
\sum_{\la\in R,i,k}\frac{(Y_\la^i)^2}{\la(A)^3}
     \la(B_0^k)\,\partial_{\al_0^k}
- \sum_{\la\in R,i}\frac{Y_\la^i}{\la(A)^2}
     \sum_{\mu,j}\langle Y,\sum_k N^{ijk}_{\la\mu} E_{\la+\mu}^k\rangle
     \partial_{Y_\mu^j}
\endalign$$
and the differential equation
$$\align
\ddot A_0^k &= 
\sum_{\la\in R,i}\frac{(Y_\la^i)^2}{\la(A)^3}\la(B_0^k),\\
\dot Y_\mu^j &= -\sum_{\la\in R,i}\frac{Y_\la^i}{\la(A)^2}
     \langle Y,\sum_k N^{ijk}_{\la\mu} E_{\la+\mu}^k\rangle 
= -\sum_{\la\in R,i}\frac{Y_\la^i}{\la(A)^2}
     \langle Y,[E_\la^i,E_\mu^j]\rangle.
\endalign$$
Let us now write 
$Y=\sum_{\la\in R,i}Y_\la^i.E_\la^i =: \sum_{\la\in R}Y_\la$ where 
$Y_\la\in\frak l_\la\cap\g$, and 
$Z=\sum_{\la\in R}\frac1{\la(A)^2}Y_\la\in \g$.
Then the differential equation becomes:
$$\align
\ddot A_0^k &= \langle \ddot A,B_0^k \rangle =
\sum_{\la\in R}\frac{\|Y_\la\|^2}{\la(A)^3}\la(B_0^k)\qquad\text{ or }
\langle \ddot A,\quad \rangle =
\sum_{\la\in R}\frac{\|Y_\la\|^2}{\la(A)^3}\la\;\in \Si^* \\
\dot Y_\mu^j &= -\Bigl\langle Y,\sum_{\la\in R,i}\frac{Y_\la^i}{\la(A)^2}
     [E_\la^i,E_\mu^j]\Bigr\rangle 
     = -\langle Y,[Z,E_\mu^j] \rangle
     = -\langle [Y,Z],E_\mu^j \rangle,
\endalign$$
so that finally we have
$$
\langle \ddot A,\quad \rangle =
\sum_{\la\in R}\frac{\|Y_\la\|^2}{\la(A)^3}\la\;\in \Si^* ,\qquad
\dot Y = -[Y,Z]. \tag{\nmb:{2}}
$$
So the ballistic curve in $C$ avoids the walls whenever 
$Y_\la(t)\ne0$ for all restricted roots $\la\in R$, and just one time 
$t$. 

\subhead\nmb.{6.6}. Symplectic reduction at singular points\endsubhead
Let us now consider a singular $A\in V$ which we may move into 
$C\subset\Si\subset V$ by using some suitable $g\in G$. Then $A$ is 
contained in some intersection of walls of $C$ so $\mu(A)=0$ for 
$\mu\in R_0\subset R$. 
Given $\al\in V$ we consider $[A,\al]=Y\in\g$. 
Using the decompositions
$\al=\al_\Si+\sum_{\la\in R, i}\al_\la^iB_\la^i$ and  
$Y=Y_{\frak h}+\sum_{\la\in R,i}Y_\la^i.E_\la^i =: 
Y_{\frak h}+\sum_{\la\in R}Y_\la$ 
as in \nmb!{6.5} we see that 
$Y=[A,\al]$ implies $Y_{\frak h}=0\in \frak h$ and 
$Y_\la^i=\la(A)\al_\la^i$ for all $\la\in R$. Thus we get $Y_\mu=0$ 
for all $\mu\in R_0$. So we can follow the analysis in \nmb!{6.5} 
without change if we agree $Y_\mu=0$ means also 
$(Y^i_\mu)^2/\mu(A)=0$ and $Z_\mu=Y_\mu/\mu(A)^2=0$. 
So again the ballistic curve is described by 
the equations \therosteritem{\nmb!{6.5.2}}.
The second equation $\dot Y=[Y,Z]$ shows that $Y_\mu=0$ along the the 
whole ballistic curve for $\mu\in R_0$. So the ballistic curve is 
composed of one just like in 
\therosteritem{\nmb!{6.5.2}} inside the intersection of walls 
$\{B\in C:\mu(B)=0\text{ for all }\mu\in R_0\}$
together with a geodesic (reflected at walls) transversal to this 
intersection of walls.


\Refs

\widestnumber\key{ABC}

\ref
\key \cit0{1}
\by Abraham, R.; Marsden, J.
\book Foundations of mechanics
\bookinfo 2nd ed
\publ Addison-Wesley
\yr 1978
\endref

\ref
\key \cit0{2}
\by Aleksandrov, A.D.
\paper \"Uber eine Verallgemeinerung der Riemann'schen Geometrie
\jour Schr. Forschungsinst. Math 
\vol 1
\yr 1957
\pages 33-84
\endref

\ref
\key \cit0{3}
\by Babelon, O.; Talon, M.
\paper The symplectic structure of the spin Calogero model.
\jour Phys. Lett. A
\vol 236, 5-6
\yr 1997
\pages 462--468
\finalinfo arXiv:q-alg/9707011
\endref

\ref
\key \cit0{4}
\by Babelon, O.; Talon, M.
\paper The symplectic structure of rational Lax pair 
\jour Phys. Lett. A
\vol 257, 3-4
\yr 1999
\pages 139--144
\endref

\ref
\key \cit0{5}
\by Bierstone, E
\paper Lifting isotopies from orbit spaces
\jour Topology
\vol 14
\yr 1975
\pages 245--252
\endref

\ref 
\key \cit0{6}
\by Burago, Yu.; Gromov, M.; Perel'man, G.
\paper A. D. Alexandrov spaces with curvature bounded below
\jour Russ. Math. Surv.
\vol 47
\yr 1992
\pages 1-58 
\endref

\ref
\key \cit0{7}
\by Chiang, Yuan-Jen
\paper Spectral geometry of $V$-manifolds and its application to 
harmonic maps
\inbook Differential geometry. Part 1: Partial differential equations 
on manifolds 
\bookinfo Proc. Symp. Pure Math. 54, Part 1
\yr 1993
\pages 93-99
\endref

\ref
\key \cit0{8}
\by Dadok, J.
\paper Polar coordinates induced by actions of compact Lie groups
\jour TAMS
\vol 288
\yr 1985
\pages 125--137
\endref

\ref
\key \cit0{9}
\by Gromov, Misha
\book Metric Structures for Riemannian and Non-Riemannian Spaces
\bookinfo Progress in Mathematics 152
\publ Birkh\"auser
\publaddr Boston, Basel, Berlin 
\yr 1999
\endref

\ref
\key \cit0{10}
\by Kazhdan, D.; Kostant, B.; Sternberg, S.
\paper Hamiltonian group actions and dynamical systems of Calogero 
type
\jour Comm. Pure Appl. Math.
\vol 31
\yr 1978
\pages 481--501
\endref

\ref
\key \cit0{11}
\by Krichever, I.M.; Babelon, O.; Biley, E.; Talon, M.
\paper  Spin generalization of the Calogero-Moser system and the 
matrix KP equation 
\inbook Topics in topology 
and mathematical physics
\bookinfo Amer. Math. Soc. Transl. Ser. 2, 170
\publ Amer. Math. Soc.
\publaddr Providence
\yr 1995
\pages 83--119
\finalinfo arXiv:hep-th/9411160
\endref

\ref
\key \cit0{12}
\by Lerman, E.; Montgomery, R.; Sjamaar, R.  
\paper Examples of singular reduction
\inbook Symplectic geometry  
\bookinfo London Math. Soc. Lecture Note Ser., 192
\publ Cambridge Univ. Press 
\publaddr Cambridge
\yr 1993
\pages 127--155
\finalinfo 95h:58054 
\endref

\ref
\key \cit0{13}
\by \L ojasiewicz, S
\paper Triangulation of semi-analytic sets 
\jour Ann. Sc. Norm. Sup. Pisa, ser. III
\vol 18
\yr 1964
\pages 449--473
\endref

\ref
\key \cit0{14}
\by Mather, J.N.
\paper Stratifications and mappings
\inbook Dynamical systems
\bookinfo New York
\publ Academic Press
\yr 1973
\pages 195-232
\endref

\ref
\key \cit0{15}
\by Mather, J.N.
\paper Differentiable invariants
\jour Topology
\vol 16
\yr 1977
\pages 145--155
\endref

\ref
\key \cit0{16}
\by Michor, Peter W.
\book Isometric actions of Lie groups and invariants
\publ Lecture course at the University of Vienna
\yr 1996/97
\endref

\ref 
\key \cit0{17} 
\by Michor, Peter W. 
\paper The Jacobi Flow 
\jour Rend. Sem. Mat. Univ. Pol. Torino
\vol 53, 3 
\yr 1996
\pages 365-372
\endref 

\ref
\key \cit0{18}
\by Moser, J.
\paper Three integrable Hamiltonian systems connected withb 
isospectral deformations
\jour Adv. Math.
\vol 16
\yr 1975
\pages 197--220
\endref

\ref
\key \cit0{19}
\by Olshanetskii, M.; Perelomov, A.
\paper Geodesic flows on symmetric spaces of zero curvature and 
explicit solution of the generalized Calogero model for the classical 
case
\jour Funct. Anal. Appl.
\vol 10, 3
\yr 1977
\pages 237--239
\endref

\ref
\key \cit0{20}
\by Palais, R.
\paper On the existence of slices for actions of non-compact Lie groups
\jour Ann. of Math. (2)
\vol 73
\yr 1961
\pages 295--323
\endref

\ref
\key \cit0{21}
\by Palais, R\. S\.; Terng, C\. L\.
\book Critical point theory and submanifold geometry
\bookinfo Lecture Notes in Mathematics 1353
\publ Springer-Verlag
\publaddr Berlin
\yr 1988
\endref

\ref
\key \cit0{22}
\by Perelomov, A.
\book Integrable systems of classical mechanics and Lie algebras
\publ Birkh\"auser
\publaddr Basel
\yr 1990
\endref

\ref
\key \cit0{23}
\by Procesi, C.; Schwarz, G.
\paper Inequalities defining orbit spaces
\jour Invent. Math.
\vol 81
\yr 1985
\pages 539--554
\endref

\ref
\key \cit0{24}
\by Satake
\paper Gauss-Bonnet theorem for $V$-manifolds
\jour J. Math. Soc. Japan
\vol 9
\yr 1957
\pages 464--492
\endref

\ref
\key \cit0{25}
\by Schwarz, G\. W\.
\paper Smooth functions invariant under the action of a compact Lie group
\jour Topology 
\vol 14
\yr 1975
\pages 63--68
\endref

\ref
\key \cit0{26}
\by Schwarz, G\. W\.
\paper Lifting smooth homotopies of orbit spaces
\jour Publ. Math. IHES
\vol 51
\yr 1980
\pages 37--136
\endref

\ref
\key \cit0{27}
\by Sjamaar, R.; Lerman, E.
\paper Stratified symplectic spaces and reduction
\jour Ann. Math.
\vol 134
\yr 1991
\pages 375--422
\endref

\ref
\key \cit0{28}
\by Whitney, H.
\paper Elementary structure of real algebraic varieties
\jour Ann. Math. 
\vol 66
\yr 1957
\pages 545--556
\endref


\endRefs
\enddocument